\documentclass[12pt]{amsart}

\usepackage{amsthm, amssymb, amsfonts}

\usepackage{mathtext}
\usepackage[cp1251]{inputenc}

\usepackage{url}

\newcommand\R{\mathbb R}

\newcommand{\Rn}{{\mathbb R}^n}
\newcommand{\Cn}{{\mathbb C}^n}

\begin{document}

\title{\bf {A generalization of the Picard theorem}}

\author{V.A.~Zorich}

\date{05.08.2021}

\begin{abstract}
We recall the notions of conformal and quasiconformal mappings \textit{in the sense of Gromov}, extending the classical notions of conformal and quasiconformal mappings, and prove the following theorem.
{\em If the mapping $ F: \R^{n} \to \R^{2} $, where $ n \geq 2 $, quasiconformal in the sense of Gromov, omits more than one value on the plane $ \R^{2} $, then it is a constant mapping.}
\end{abstract}



\maketitle
\markboth{Vladimir Zorich}{Generalized Picard theorem}

{\bf 1.~Introduction.}
~M.Gromov extended the concepts of conformal and quasiconformal mapping
to the mappings acting between the manifolds of different dimensions. For instance, any entire holomorphic function $ f: \Cn \to {\mathbb C}$ defines a mapping conformal in the sense of Gromov.

A mapping $ F: \R^{m} \to \R^{n} $ ~($ {m} \geq {n} $)
is {\em conformal} in the sense of Gromov
if it transformes any infinitesimally small ball of domain into an infinitesimally small ball of the image [1].

A mapping is {\em quasiconformal} in the sense of Gromov if the image of any
infinitesimal ball is an infinitesimal ellipsoid of eccentricity
bounded from above by a constant independent of the point.

\vskip2mm

In connection with such an extension of the concepts of conformality and quasiconformality of a mapping, Gromov naturally raises the question of what facts of the classical complex analysis apply to these mappings. For example, is it true that
{\em if the mapping $ F: \R^{m} \to \R^{n} $ is conformal and bounded, then for
$ m \geq n \geq 2 $
it is a constant ?}

In the paper [2] we discussed this issue, and later showed in~[3] that, as Gromov probably expected, the generalization of the classical Liouville theorem, stated above as a question, is indeed valid. We have shown that such a Liouville constant theorem holds not only for conformal but also for quasiconformal mappings.

\vskip2mm

Below we prove that the corresponding generalization of the classical Picard theorem also takes place.

It is curious that the development of the theory of quasiconformal mappings itself once began with Gr\"otzsch's observation [4] that the seemingly specifically holomorphic Picard theorem turned out to be valid for a much wider class of mappings which Ahlfors later called {\em quasiconformal mappings}.

\vskip2mm

{\bf 2.~The Picard theorem.}~
In his famous book
``Mathematical Miscellany''~
[5], discussing the question of whether there are great theorems in mathematics with very short proofs, Littlewood cites the Picard theorem as an example. I quote Littlewood's words full of humor and professionalism.

\textit{
``The question recently arose in conversation whether a dissertation of 2 lines could deserve and get
a Fellowship. I had answered this for myself long before; in mathematics the answer is yes.
}

\textit{
... With Picard's Theorem it could be literally 2, one of statement, one of proof.}

\vskip1mm

\textit{Theorem.
An integral function never 0 or 1 is a constant.}
\vskip1mm

\textit{
Proof.
$\exp (i \Omega (f (z)))$ is a bounded integral function.}

\vskip2mm

\textit{
...
With Picard the situation is clear enough today (innumerable papers have resulted from it). But I can imagine a referee's report: `Exceedingly striking and a most original idea. But, brilliant as it undoubtedly is, it seems more odd than important; an isolated result, unrelated to anything else, and not likely to lead anywhere.''}

\vskip2mm
The given single  line of the complete proof of the Picard theorem in terms of conformal mappings, modular function, and the monodromy theorem, of course, assumes that the reader is sufficiently advanced in complex analysis. These classical arguments reduce the Picard theorem to the Liouville constant theo\-rem. They are also valid in the considered case of mappings conformal and quasiconformal in the sense of Gromov if the aforementioned Liouville theorem has already been established for them.
Thus, the following theorem holds.

\vskip2mm

{\bf Theorem.}
{\em If the mapping $ F: \R^{n} \to \R^{2} $, ($ n \geq 2 $) quasiconformal in the sense of Gromov omits more than one value on the plane $ \R^{2} $, then it is a constant mapping. }

\vskip2mm

{\bf 3.~Remark.}~
It is well known that
Picard theorem is valid for entire holomorphic functions
$ f: \Cn \to {\mathbb C} $.
We have shown above that the theorem is also valid for mappings
$ f:  \Rn \to {\mathbb R}^2 $ conformal or quasiconformal in the sense of Gromov if $ n \geq 2 $ .

It is also well known that for holomorphic mappings
$ f: \Cn \to {\mathbb C}^m $
Liouville constant theorem is valid, while
Picard theorem for $ m> 1 $, generally speaking, does not hold. The set $ {\mathbb C}^m \setminus f (\Cn) $ can even
contain entire neighborhoods of some points (recall Fatou's example for the case $ n = m = 2 $).

At the same time, for quasiconformal mappings $ f: \Rn \to {\mathbb R}^n $ and $ n > 1 $
Picard theorem is valid in the following form:
{\em if quasiconformal mapping (generally speaking,
not univalent)
is not a constant  mapping, then the set
$ \Rn \setminus f (\Rn) $ is finite, and the number of its points does not exceed some value $ P (n, k_f) $, depending only on the dimension of the space and on the quasiconformality coefficient of the mapping} [6], [7].

One might expect that
Picard theorem of this kind is also valid for mappings quasiconformal in the sense of Gromov.

\vskip3mm

\centerline {\bf References}
\vskip2mm

[{\bf 1}] M.L.~Gromov, “Colourful categories”, Russian Math. Surveys, 70:4 (2015), 591–655.
{\footnote{For the information of the reader we also note that on the very first page of the article [1], in a footnote, the author gives the address where a much more complete text of his work
is presented. In particular, there the reader can find an extended interpretation of conformity and quasiconformality, as well as the question about the Liouville theorem for such mappings. The text is now available at the site  \small{\url {https://www.ihes.fr/~gromov/wp-content/uploads/2018/08/problems-sept2014-copy.pdf}}}.

[{\bf 2}] V.A.~Zorich, “On a question of Gromov concerning the generalized Liouville theorem”,
Russian Math. Surveys, 74:1 (2019), 175–177.

[{\bf 3}] V.A.~Zorich, “Conformality in the sense of Gromov and a generalized Liouville theorem”,
http://arxiv.org/abs/2108.00945

[{\bf 4}] H.~Gr{\" o}tzsch,
{\" U}ber die Verzerrung bei schlichten nichtkonformen Abbil\-dun\-gen
und {\" u}ber eine damit zusammenh{\" a}ngende Erweiterung des Picardschen Satzes.
Ber. Verh. S{\" a}chs.
Akad. Wiss. Leipzig. 1928. V. 80. P. 503-507.

[{\bf 5}] J.E.~Littlewood, “A Mathematical Miscellany”. London, 1957.

[{\bf 6}] S.~Rickman, “Quasiregular Mappings”,
A Series of Modern Surveys in Mathematics, {\bf 26},
Springer-Verlag Berlin Heidelberg, 1993.

[{\bf 7}] A.~Eremenko, “Value distribution and potential theory”,
Proceedings of the ICM, Beijing 2002, vol. 2, 681--690.
Higher Ed. Press, Beijing, 2002.
%
%
%

\vskip3mm

{\small {\bf V.A.~Zorich}

Lomonosov Moscow State University

\vskip1mm
E-mail: {\bf vzor@mccme.ru}

\end{document}